\documentclass[11pt]{article}
\usepackage{authblk}
\usepackage{amssymb}
\usepackage[all]{xy}
\usepackage{mathrsfs}
\usepackage{amsmath,amssymb,amsthm}
\usepackage{extarrows}
\usepackage{mathrsfs}
\newtheorem{theorem}{Theorem}[section]
\newtheorem{lemma}[theorem]{Lemma}
\newtheorem{proposition}[theorem]{Proposition}
\newtheorem{definition}[theorem]{Definition}

\def\to{\rightarrow}
\def\fA{\mathfrak{A}}
\def\fB{\mathfrak{B}}
\def\fC{\mathfrak{C}}
\def\fQ{\mathfrak{Q}}
\def\mA{\mathcal{A}}

\def\mC{\mathcal{C}}
\def\bC{\mathbf{C}}
\def\ot{\otimes}
\def\P{\Phi}
\def\r{\rho}
\def\m{\mu}
\def\n{\nu}
\date{}
\begin{document}
\title{Quantum metric spaces of quantum maps}
\author{Maysam Maysami Sadr\thanks{sadr@iasbs.ac.ir}}
\affil{Department of Mathematics, Institute for Advanced Studies in Basic Sciences, Zanjan, Iran}
\maketitle
\begin{abstract}
We show that any quantum family of quantum maps from a non commutative space to a compact quantum metric
space has a canonical quantum pseudo-metric structure.


\end{abstract}
\maketitle
\section{Introduction}
One of the basic ideas of {\it Noncommutative Geometry} is that any unital C*-algebra $A$ can be considered as the
algebra of {\it continuous functions} on a  (symbolic) {\it compact quantum (noncommutative) space} $\fQ A$. From this point of view,
any unital *-homomorphism $\P:B\to A$ between unital C*-algebras can be interpreted as a {\it quantum map} $\fQ\P$ from $\fQ A$ into $\fQ B$.
There are many notions in Topology and Geometry that can be translate into NC language.
The notion of {\it quantum family of (quantum) maps}, defined by Woronowicz \cite{W1} and So{\l}tan \cite{Soltan1} (see also \cite{Sadr3,Sadr2,Sadr1}),
conclude from the following fact:
``Every map $f$ from $X$ to the set of all maps from $Y$ to $Z$ (or in other word, any  family of maps
from $Y$ to $Z$  parameterized by $f$ with  parameters $x$ in $X$) can be considered as a map
$\tilde f:X\times Y\to Z$ defined by $\tilde f(x,y)=f(x)(y)$.'' A translation of this to noncommutative language is as follows.
\begin{definition}
(\cite{Sadr3,Sadr2,Sadr1,Soltan1,W1} Let $B,C$ be unital C*-algebras. A quantum family of morphisms from $B$ to $C$ (or, a quantum family of maps
from $\fQ C$ to $\fQ B$) is a pair $(A,\P)$ consisting of a unital C*-algebra $A$ and a unital *-homomorphism
$\P:B\to C\ot A$, where $\ot$ denotes the spatial tensor product of C*-algebras.
\end{definition}
Another concept that can be translate from Geometry into NC Geometry, is {\it distance} or {\it metric}.
Marc Rieffel, by using the notion of {\it order unite spaces}, has developed the notion of {\it quantum metric
space} in a series of papers \cite{R1,R2,R3,R4,R6}. For two other different notions of quantum metric see \cite{KuperbergWeaver1,Sadr4,Sadr6,Sadr5}.
Here, we deals with special examples of Rieffel's quantum metric spaces, stated in the C*-algebraic formalism.
The aim of this note is to show that any quantum family of maps from a quantum space to a compact quantum
metric space has a canonical quantum pseudo-metric structure. We are motivated by the following trivial fact:
Let $(Z,d)$ be a metric space and $f:X\times Y\to Z$ be a family of maps from $Y$ to $Z$, then $X$
has a pseudo-metric $\r$ defined by $$\r(x,x')=\sup_{y\in Y}d(f(x,y),f(x',y)).$$
In Section \ref{s2} we introduce the notion of {\it compact quantum pseudo-metric space}. In Section \ref{s3} we define a natural
compact quantum pseudo-metric space structure on any quantum family of maps from a quantum space to a compact quantum metric space.
In Section \ref{s4} we examine our definition in the classical case.
\section{Compact quantum pseudo-metric spaces}\label{s2}
By a pseudo-metric $d$ on a set $X$ we mean a positive valued function on $X\times X$ which is symmetric, satisfies triangle inequality,
and $d(x,x)=0$ for every $x\in X$. For any topological space $X$ with topology $\tau$ (resp. pseudo-metric space $(X,d)$) $\bC(X,\tau)$ (resp.
$\bC(X,d)$) denotes the C*-algebra of all continuous bounded complex valued maps on $X$ with the uniform norm.
For a pseudo-metric $d$, $\tau_d$ denotes the topology induced by $d$. Let $(X,d)$ be a
pseudo-metric space. For every $f\in\bC(X,d)$, the Lipschitz semi norm $\|f\|_d$ is defined by
$$\|f\|_d=\sup\{\frac{|f(x)-f(x')|}{d(x,x')}: x,x'\in X, d(x,x')\neq0\}.$$
Also, the Lipschitz algebra of $(X,d)$ is defined by,
$${\bf Lip}(X,d)=\{f\in\bC(X,d):\|f\|_d<\infty\}.$$
We need the following simple lemma.
\begin{lemma}\label{l1}
Let $(X,d)$ be a pseudo-metric space and $a$ be a complex valued map on $X$. Then $a\in{\bf Lip}(X,d)$ and
$\|a\|_d\leq1$ if and only if $|a(x)-a(x')|\leq d(x,x')$ for every $x,x'\in X$. In particular, if $b\in\bC(X,d)$, then $\|b\|_d=0$ if and only if $b$ is a constant map.
\end{lemma}
\begin{proof}
Let $a\in{\bf Lip}(X,d)$ and $\|a\|_d\leq1$. Suppose that $x,x'\in X$. If $d(x,x')=0$, then $a(x)=a(x')$, since
$a$ is continuous with $\tau_d$. If $d(x,x')\neq0$, then $1\geq\|a\|_d\geq\frac{|a(x)-a(x')|}{d(x,x')}$, and thus
$|a(x)-a(x')|\leq d(x,x')$. The other direction is trivial.
\end{proof}
For any C*-algebra $\fA$, $S(\fA)$ denotes the state space of $\fA$ with w* topology.
If $\fA$ is  unital, $1_\fA$ denotes the unit element of $\fA$.\\
Let $\mA$ be a self adjoint linear subspace of the C*-algebra $\fA$, and let $L:\mA\to[0,\infty)$ be a
semi norm on $\mA$. Connes has pointed out \cite{C1}, \cite{C2}, that one can define a pseudo-metric
$\r_L$ on $S(\fA)$ by
\begin{equation}\label{e1}
\r_L(\m,\n)=\sup\{|\m(a)-\n(a)|: a\in\mA,L(a)\leq1\}\quad\quad(\m,\n\in S(\fA)).
\end{equation}
Note that $\r_L$ can take values $+\infty$ and $0$ for different states of $\fA$. Conversely, let $d$ be a pseudo-metric
on $S(A)$ (such that the topology induced by $d$ on $S(\fA)$ is not necessarily w* topology). Define
a semi norm $L_d:\fA\to[0,+\infty]$ by
$$L_d(a)=\sup\{\frac{|\m(a)-\n(a)|}{d(\m,\n)}:\m,\n\in S(\fA), d(\m,\n)\neq0\}\quad\quad(a\in\fA).$$
Note that $L_d(a)=L_d(a^*)$ for every $a\in\fA$.\\
Let $(X,d)$ be a compact metric space. Consider the Lipschitz semi norm
$$\|\cdot\|_d:{\bf Lip}(X,d)\subset\bC(X,d)\to[0,+\infty).$$
Then it is easily checked that the semi norm $\r_{\|\cdot\|_d}$ on the state space of $\bC(X,d)$ is a metric,
called Monge-Kantorovich metric \cite{Ra}.
It is well known that the topology induced by $\r_{\|\cdot\|_d}$, is the w* topology, and for every $x,y\in X$,
$d(x,y)=\r_{\|\cdot\|_d}(\delta_x,\delta_y)$, where $\delta:X\to\bC(X,d)^*$ is the point mass measure map.
\begin{proposition}\label{p1}
Let $(X,\tau)$ be a compact Hausdorff space and $d$ be a pseudo-metric on $X$ such that the topology induced
by $d$ on $X$ is weaker than $\tau$, i.e. $\tau_d\subset\tau$. Consider the Lipschitz semi norm
$\|\cdot\|_d:{\bf Lip}(X,d)\subset\bC(X,\tau)\to[0,+\infty)$ and let $\r=\r_{\|\cdot\|_d}$. Then the following are
satisfied.
\begin{enumerate}
\item [i)] $d(x,y)=\r(\delta_x,\delta_y)$, for every $x,y\in X$.
\item [ii)] $L_\r=\|\cdot\|_d$ on $\bC(X,d)\subset\bC(X,\tau)$.
\item [iii)] Let $a\in\bC(X,\tau)$, then $a\in\bC(X,d)$ if and only if the map
$\n\longmapsto\n(a)$ on $S(\bC(X,\tau))$ is continuous with $\r$.
\item [iv)] the topology induced by $\r$ on $S(\bC(X,\tau))$ is weaker than the w* topology.
\end{enumerate}
\end{proposition}
\begin{proof}
i) Let $x,y$ be in $X$. Suppose that $a\in{\bf Lip}(X,d)$ and $\|a\|_d\leq1$. Then by Lemma \ref{l1},
$|\delta_x(a)-\delta_y(a)|=|a(x)-a(y)|\leq d(x,y)$, and thus by definition of $\r$, we have
$\r(\delta_x,\delta_y)\leq d(x,y)$. Conversely, let $a_x\in\bC(X,d)$ be defined by $a_x(z)=d(x,z)$ ($z\in X$);
then for every $x',y'\in X$, $|a_x(x')-a_x(y')|=|d(x,x')-d(x,y')|\leq d(x',y')$, and thus by lemma \ref{l1},
$a\in{\bf Lip}(X,d)$ and $\|a\|_d\leq1$. Now, we have
$$\r(\delta_x,\delta_y)\geq|\delta_x(a_x)-\delta_y(a_x)|=|a_x(x)-a_x(y)|=d(x,y).$$

ii) By i) and definitions of $L_\r$ and $\|\cdot\|_d$, it is clear that $\|\cdot\|_d\leq L_\r$ on $\bC(X,\tau)$.\\
Let $a\in\bC(X,d)$. If $\|a\|_d=0$, then by Lemma \ref{l1}, $a$ is a constant map and thus $L_\r(a)=0$.
If $\|a\|_d=\infty$ then $L_\r(a)=\infty$ since $\|a\|_d\leq L_\r(a)$. Thus suppose that $0<\|a\|<\infty$.
Then for every $\m,\n\in S(\bC(X,\tau))$, we have
$$\r(\m,\n)\geq|\m(\frac{a}{\|a\|_d})-\n(\frac{a}{\|a\|_d})|=\frac{|\m(a)-\n(a)|}{\|a\|_d}$$
and thus if $\r(\m,\n)\neq0$ then $\|a\|_d\geq\frac{|\m(a)-\n(a)|}{\r(\m,\n)}$. Therefore,
$$\|a\|_d\geq\sup\{\frac{|\m(a)-\n(a)|}{\r(\m,\n)}:\hspace{2mm}\m,\n\in S(\bC(X,\tau)),\r(\m,\n)\neq0\}=L_\r(a).$$

iii) The `if' part is an immediate consequence of i). For the other direction, we need some notations:
Let $\sim$ be the equivalence relation on $X$ defined by $x\sim x'\Leftrightarrow d(x,x')=0$.
Let $Y=X/\sim$ and let $\hat{}:X\to Y$ be the canonical projection.  Then $\hat d$, defined by
$\hat d(\hat x_1,\hat x_2)=d(x_1,x_2)$, is a well defined metric on $Y$, and $\hat{}$ is an isometry between
$(X,d)$ and $(Y,\hat d)$. Thus the C*-algebras $\bC(X,d)$ and $\bC(Y,\hat d)$, and the Lipschitz algebras
$({\bf Lip}(X,d),\|\cdot\|_{d})$ and $({\bf Lip}(Y,\hat d),\|\cdot\|_{\hat d})$ are isometric isomorph. In particular,
the topology induced by $\r$ on $S(\bC(X,d))$ is the w* topology, since as mentioned above the
Monge-Kantorovich metric $\r_{\|\cdot\|_{\hat d}}$ induces the w* topology on $S(\bC(Y,\hat d))$. Consider the
canonical embedding $\P:\bC(X,d)\to\bC(X,\tau)$. For every $\n,\n'\in S(\bC(X,\tau))$, $\n\circ\P$ and
$\n'\circ\P$ are in $S(\bC(X,d))$ and
\begin{equation}\label{e9}
\r(\n,\n')=\r(\n\circ\P,\n'\circ\P).
\end{equation}
Now, let $a\in\bC(X,d)$ and $\n_i\to\n$ be a convergent net in $S(\bC(X,\tau))$ with $\r$. Then
$\n_i\circ\P\to\n\circ\P$ is a convergent net in $S(\bC(X,d))$ with $\r$, and since the topology induced by $\r$
agrees with the w* topology on $S(\bC(X,d))$, we have
$$\n_i(a)=\n_i\circ\P(a)\to\n\circ\P(a)=\n(a).$$
Thus we get the desired result.

iv) Let $\n_i\to\n$ be a convergent net in $S(\bC(X,\tau))$ with w* topology. Thus as in the proof of iii),
$\n_i\circ\P\to\n\circ\P$ with $\r$, and by (\ref{e9}),  $\n_i\to\n$ in $S(\bC(X,\tau))$ with the topology
induced by $\r$. This completes the proof of iv).
\end{proof}
\begin{definition}
By a compact quantum pseudo-metric space (QSM space, for short) we mean a triple $(\fA,\mA,L)$,
where $\fA$ is a unital C*-algebra, $\mA$ is a
self adjoint linear subspace of $\fA$ with $1_\fA\in\mA$, and $L:\mA\to[0,+\infty)$ is a semi norm such that
\begin{enumerate}
\item [(a)] $L(a)=L(a^*)$ for every $a\in\mA$,
\item [(b)] for every $a\in\mA$, $L(a)=0$ if and only if $a\in\mathbb{C}1_\fA$, and
\item [(c)] the topology induced by the pseudo-metric $\r_L$ on $S(\fA)$ is weaker than the w*
topology.
\end{enumerate}
\end{definition}
As an immediate corollary of the definition, for any compact quantum pseudo-metric space $(\fA,\mA,L)$, the
topology induced by $\r_L$ on $S(\fA)$ is compact and in particular the diameter of $S(\fA)$ under $\r_L$ is finite.
\begin{proposition}\label{p3}
Let $(\fA,\mA,L)$ be a QSM space. Then, for every $a\in\mA$, the map $\m\longmapsto\m(a)$ on $S(\fA)$ is
continuous with topology induced by $\r_L$.
\end{proposition}
\begin{proof}
Straightforward.
\end{proof}
\begin{definition}
A QSM space $(\fA,\mA,L)$ is called a compact quantum metric space (QM space, for short) if $\mA$ is a dense
subspace of $\fA$.
\end{definition}
Let $(\fA,\mA,L)$ be a QM space and $\m,\n$ be two different states of $\fA$. Then since
$\mA$ is dense in $\fA$, there is $a\in\mA$ such that $\m(a)\neq\n(a)$. Thus (by (\ref{e1})) $\r_L$ is a
metric on $S(\fA)$. It is an elementary result in Topology that any Hausdorff topology $\tau$ weaker than
a compact Hausdorrf topology $\tau'$ on a set $X$, is equal to the same topology $\tau'$. Using this, we conclude
that the topology induced by $\r_L$ on $S(\fA)$ is the w* topology.

{\bf Example 1.} Let $(X,d)$ be a compact metric space. Then
$$(\bC(X,d),{\bf Lip}(X,d),\|\cdot\|_d)$$ is a compact
quantum metric space.

{\bf Example 2.} Let $(X,\tau)$ be a compact Hausdorff space and let $d$ be a pseudo-metric on $X$ such that
$\tau_d\subset\tau$. Then Proposition \ref{p1} and Lemma \ref{l1}, show
$$(\bC(X,\tau),{\bf Lip}(X,d),\|\cdot\|_d)$$
is a compact quantum pseudo-metric space.

{\bf Remark.} Let $(\fA,\mA,L)$ be a QM space and $A\subset\mA$ be the linear subspace of
all self-adjoint elements of $\mA$. Then $A$ is an order unite space and $(A,L|_A)$ is a compact quantum
metric space in the sense of Riffel's definition \cite{R3}.
\begin{lemma}\label{l4}
Let $\fA$ be a C*-algebra with the C*-norm $\|\cdot\|$, $\mA$ be a self adjoint linear subspace of $\fA$ containing
$1_\fA$ and $L:\mA\to[0,+\infty)$ be a semi norm such that for every $a\in\mA$, $L(a)=0$ if and only if
$a\in\mathbb{C}1_\fA$. Let $\tilde L$ and $\|\cdot\|\tilde{}$ denote the quotient norm of $L$ and $\|\cdot\|$ on
$\frac{\mA}{\mathbb{C}1_\fA}$ and $\frac{\fA}{\mathbb{C}1_\fA}$, respectively. Suppose that the image of
$\{a\in\mA:L(a)\leq1\}$ in $\frac{\fA}{\mathbb{C}1_\fA}$ is totally bounded for $\|\cdot\|\tilde{}$. Then the
topology induced by $\r_L$ on $S(\fA)$ is weaker than the w* topology.
\end{lemma}
\begin{proof}
See Theorem 1.8 of \cite{R1}.
\end{proof}
{\bf Example 3.} Let $\fA$ be a finite dimensional C*-algebra and $N$ be a Banach space norm on $\fA$ such that
$N(a)=N(a^*)$ for every $a\in\fA$. Let the semi norm $N_0:\fA\to[0,\infty)$ be defined by
$$N_0=\inf\{N(a+\lambda1_\fA):\hspace{3mm}\lambda\in\mathbb{C}\}.$$
Since $\fA$ is finite dimensional, the C*-norm of $\fA$ and $N$ are equivalent. Thus the image $K$ of
$\{a\in\fA:N_0(a)\leq1\}$ is closed and bounded in $\frac{\fA}{\mathbb{C}1_\fA}$. Again, since $\fA$ is finite
dimensional, $K$ is compact and thus totally bounded for the quotient norm of the C*-norm. Thus by Lemma \ref{l4},
$(\fA,\fA,N_0)$ is a QM space.

{\bf Example 4.} Let $G$ be a compact Hausdorff group with identity element $e$. Let $\ell$ be a length function
on $G$, i.e. $\ell$ is a continuous non negative real valued function on $G$ such that
\begin{enumerate}
\item [(i)] $\ell(gg')\leq\ell(g)+\ell(g')$, for every $g,g'\in G$,
\item [(ii)] $\ell(g)=\ell(g^{-1})$ for every $g\in G$, and
\item [(iii)] $\ell(g)=0$ if and only if $g=e$.
\end{enumerate}
Let $\fA$ be a unital C*-algebra with a strongly continuous action $\cdot:G\times\fA\to\fA$ of $G$ by automorphisms
of $\fA$, i.e.
\begin{enumerate}
\item [(a)] for every $g\in G$ the map $a\longmapsto g\cdot a$ is a *-automorphism of $\fA$,
\item [(b)] $e\cdot a=a$ for every $a\in\fA$ ,
\item [(c)] $g\cdot(g'\cdot a)=(gg')\cdot a$, for every $g,g'\in G, a\in A$, and
\item [(d)] if $g_i\to g$ is a convergent net in $G$ and $a\in\fA$, then $g_i\cdot a\to g\cdot a$ with the C*-norm
of $\fA$.
\end{enumerate}
Define a semi norm $L$ on $\fA$ by
$$L(a)=\sup\{\frac{\|g\cdot a-a\|}{\ell(g)}:\hspace{2mm}g\in G,g\neq e\}\hspace{10mm}(a\in\fA).$$
Let $\mA=\{a\in\fA:L(a)<+\infty\}$. Then by Proposition 2.2 of \cite{R1}, $\mA$ is a dense *-subalgebra of $\fA$.
Now, suppose that the action of $G$ is ergodic, i.e.  if $a\in\fA$ and for every $g\in G$, $g\cdot a=a$, then
$a\in\mathbb{C}1_\fA$. Then it is trivial that $L(a)=0$ if and only if $a\in\mathbb{C}1_\fA$. Rieffel has proved
\cite[Theorem 2.3]{R1}, that the topology induced by $\r_L$ on $S(\fA)$ agrees with the w* topology. Thus
$(\fA,\mA,L)$ is a QM space.

For some other examples that completely match our notion of QM space, see \cite{R1}.
As we will see in the next section, using quantum family of
morphisms we can construct many QSM spaces from a QSM space.
\section{The main definition}\label{s3}
We need the following simple topological lemma.
\begin{lemma}\label{l2}
Let  $Y$ be a compact space, $X$ be an arbitrary space and $(Z,\r)$ be a pseudo-metric space. Also, let $\bC(Y,Z)$
be the space of all continuous maps from $Y$ to $Z$, with the pseudo-metric $\hat\r$ defined by
$$\hat\r(f,g)=\sup\{\r(f(y),g(y)):\hspace{5mm}y\in Y\}\hspace{10mm}(f,g\in\bC(Y,Z)).$$
Suppose that $F:Y\times X\to Z$ is a continuous map. Then the map $\tilde F:X\to\bC(Y,Z)$, defined by
$\tilde F(x)(y)=F(y,x)$ is continuous.
\end{lemma}
\begin{proof}
Let $x_0\in X$ and $\epsilon>0$ be arbitrary. Since $F$ is continuous, for every $y\in Y$, there are open
sets $U_y, V_y$ in $X$ and $Y$ respectively, such that $(y,x_0)\in V_y\times U_y$ and
$\r(F(y,x_0),F(y',x))<\epsilon/2$ for every $(y',x)\in V_y\times U_y$. Since $Y$ is compact, there are
$y_1,\cdots,y_n\in Y$ such that $Y=\cup_{i=1}^nV_{y_i}$. Let $W$ be the open set $\cap_{i=1}^nU_{y_i}$.
Let $x\in W$ and $y\in Y$ be arbitrary. Then for some $i$ ($i=1,\cdots,n$), $y$ belongs to $V_{y_i}$ and we have,
$$\r(F(y,x),F(y,x_0))\leq\r(F(y,x),F(y_i,x_0))+\r(F(y_i,x_0),F(y,x_0))<\epsilon.$$
Thus we have $\hat\r(\tilde F(x),\tilde F(x_0))<\epsilon$ for every $x\in W$. The proof is complete.
\end{proof}
Let $(\fA,\mA,L)$ be a QSM space, $\fB$ be a unital C*-algebra, and $(\fC,\P)$ be a
quantum family of morphisms from $\fA$ to $\fB$, $\P:\fA\to\fB\ot\fC$.\\
Let $d$ be a pseudo-metric on $S(\fC)$, defined  by
$$d(\n,\n')=\sup\{\r_L((\m\ot\n)\P,(\m\ot\n')\P):\hspace{2mm}\m\in S(\fB)\}\hspace{10mm}(\n,\n'\in S(\fC)).$$
\begin{proposition}\label{p2}
With the above assumptions, let $\mC$ be the linear space of all $c\in\fC$ such that the map
$\n\longmapsto\n(c)$ on $S(\fC)$ is continuous with the topology induced by $d$, and $L_d(c)<\infty$.
Then the following are satisfied.
\begin{enumerate}
\item [i)] $\mC$ is a self adjoint linear subspace of $\fC$ and $1_\fC\in\mC$.
\item [ii)] For every $c\in\mC$, $L_d(c)=0$ if and only if $c\in\mathbb{C}1_\fC$.
\item [iii)] The topology induced by $d$ on $S(\fC)$ is weaker than the w* topology.
\item [iv)] With the restriction of the domain of $L_d$ to $\mC$, $\r_{L_d}\leq d$.
\item [v)] The topology induced by $\r_{L_d}$ on $S(\fC)$ is weaker than the w* topology.
\end{enumerate}
\end{proposition}
\begin{proof}
i) is easily checked.

ii) Let $c$ be in $\mC$ and $L_d(c)=0$. By Lemma \ref{l1}, the map $\n\longmapsto\n(c)$ on $S(\fC)$ is constant,
and thus $c\in\mathbb{C}1_\fC$.

iii) Apply Lemma \ref{l2}, with $X=S(\fC)$, $Y=S(\fB)$, $Z=S(\fA)$, $\r=\r_L$ and $F:Y\times X\to Z$ defined by
$$F(\m,\n)=(\m\ot\n)\P\hspace{10mm}(\m\in Y,\n\in X).$$
We get $\tilde F:X\to \bC(Y,Z)$ is continuous with the metric $\hat\r$ on $\bC(Y,Z)$. On the other hand, for every
$\n,\n'$ we have $d(\n,\n')=\hat\r(\tilde F(\n),\tilde F(\n'))$. Thus, if $\n_i\to\n$ is a convergent net in $X$
with w* topology, then
$$d(\n_i,\n)=\hat\r(\tilde F(\n_i),\tilde F(\n))\to0.$$
This implies that the topology induced by $d$ is weaker than the w* topology.

iv) Let $\n,\n'$ be in $S(\fC)$. If $d(\n,\n')=0$ then for every $c\in\mC$, $\n(c)=\n'(c)$
(since the map $\m\longmapsto\m(c)$ is continuous with $d$) and thus by the definition of $\r_{L_d}$,
$\r_{L_d}(\n,\n')=0$. Thus suppose that $d(\n,\n')\neq0$. Let $c\in\mC$ with $L_d(c)\leq1$. Then
$1\geq L_d(c)\geq\frac{|\n(c)-\n'(c)|}{d(\n,\n')}$, and thus $|\n(c)-\n'(c)|\leq d(\n,\n')$. Therefore
$$\r_{L_d}(\n,\n')\leq d(\n,\n').$$

v) follows directly from iv) and iii).
\end{proof}
\begin{definition}
With the above assumptions, Proposition \ref{p2}, shows that
$(\fC,\mC,L_d)$ is a QSM  space that is called QSM space induced by
the QSM space $(\fA,\mA,L)$ and quantum family of maps $(\fC,\P)$.
\end{definition}
\begin{lemma}\label{l3}
With the above assumptions, let $a\in\mA$ and let $\m\in S(\fB)$. Then $c=(\m\ot id_\fC)\P(a)$ is in $\mC$, and
$L_d(c)\leq L(a)$.
\end{lemma}
\begin{proof}
We first show that $L_d(c)\leq L(a)(<\infty)$. If $L(a)=0$ then $a\in\mathbb{C}1_\fA$ and thus $c\in\mathbb{C}1_\fC$
and $L_d(c)=0$. Suppose that $L(a)\neq0$. We prove that for every $\n,\n'\in S(\fC)$ with $d(\n,\n')\neq0$,
\begin{equation}\label{e2}
\frac{|\n(c)-\n'(c)|}{d(\n,\n')}\leq L(a).
\end{equation}
Let $\n,\n'\in S(\fC)$ be such that $d(\n,\n')\neq0$. If $|\n(c)-\n'(c)|=0$, then (\ref{e2}) is satisfied.
Suppose that
$$|\n(c)-\n'(c)|=|(\m\ot\n)\P(a)-(\m\ot\n')\P(a)|\neq0.$$
By the definition of $d$, we have $d(\n,\n')\geq\r_L((\m\ot\n)\P,(\m\ot\n')\P)$. On the other hand, by the definition
of $\r_L$,
\begin{equation*}
\begin{split}
\r_L((\m\ot\n)\P,(\m\ot\n')\P)&\geq|(\m\ot\n)\P(\frac{a}{L(a)})-(\m\ot\n')\P(\frac{a}{L(a)})|\\
&=\frac{|(\m\ot\n)\P(a)-(\m\ot\n')\P(a)|}{L(a)}.
\end{split}
\end{equation*}
Thus, (\ref{e2}) is satisfied and $L_d(c)\leq L(a)$.

Now, we show that the map $\n\longmapsto\n(c)$ on $S(\fC)$ is continuous with $\tau_d$. Let $\n_n\to\n$ be a
convergent sequence in $S(\fC)$ with the metric $d$. Thus, by the definition of $d$, we have
$$\r_L((\m\ot\n_n)\P,(\m\ot\n)\P)\to0.$$
Therefore, by Proposition \ref{p3},
$$\n_n(c)=(\m\ot\n_n)\P(a)\to(\m\ot\n)\P(a)=\n(c).$$
\end{proof}
\begin{proposition}\label{p4}
With the above assumptions, suppose that $(\fA,\mA,L)$ is a QM space and the linear span of
$$G=\{(\m\ot id_\fC)\P(a):\hspace{3mm}\m\in S(\fB), a\in\fA\}$$
is dense in $\fC$ (for example $\P$ is surjective). Then $(\fC,\mC,L_d)$ is a QM space.
\end{proposition}
\begin{proof}
Since $\mA$ is dense in $\fA$ and the linear span of $G$ is dense in $\fC$, we have
$$G_0=\{(\m\ot id_\fC)\P(a):\hspace{3mm}\m\in S(\fB), a\in\mA\}$$
is dense in $\fC$. On the other hand, by Lemma \ref{l3}, $G_0\subset\mC$. Thus $\mC$ is dense in $\fC$ and
$(\fC,\mC,L_d)$ is a QM space.
\end{proof}
{\bf Example 5.} Let $\fA$ and $\fC$ be unital C*-algebras. Suppose that $\fA\ot\fC$ has a QSM structure.
Consider *-homomorphisms
$$id:\fA\ot\fC\to\fA\ot\fC\hspace{4mm}\text{and}\hspace{4mm}F:\fA\ot\fC\to\fC\ot\fA,$$
where $F$ is the {\it flip} map, i.e. $F(a\ot c)=c\ot a$  for $a\in\fA,c\in\fC$. Then
$$(\fC,id_{\fA\ot\fC})\hspace{10mm}\text{and}\hspace{10mm}(\fA,F)$$
are quantum families of morphisms. Thus $\fA$ and $\fC$ have naturally QSM structures.
Also, by Proposition \ref{p4}, if $\fA\ot\fC$ has a QM structure then so are $\fA$ and $\fC$.

{\bf Example 6.} Let $\fA$ be a unital C*-algebra and suppose that $\fA$ has a QSM structure. Let $\P:\fA\to\fB$
be a unital *-homomorphism. Then $(\fB,\P)$ can be considered as a quantum family of morphisms from $\fA$ to
$\mathbb{C}$. Thus $\fB$ naturally has a QSM structure. Also, if $\P$ is surjective and $\fA$ has a QM structure,
then by Proposition \ref{p4}, $\fB$ has a QM structure.
\section{The commutative case}\label{s4}
In this last section we study induced metric structures on ordinary families of maps.
\begin{lemma}\label{l5}
Let $(X,\tau)$ be a compact Hausdorff space and let $d$ be a pseudo-metric on $S(\bC(X,\tau))$ such that $\tau_d$
is weaker than the w* topology. Let $\mC$ be the space of all $c\in\bC(X,\tau)$ such that the map
$\n\longmapsto\n(c)$ is continuous on $S(\bC(X,\tau))$ and $L_d(c)<\infty$. Consider the semi norm
$L_d:\mC\to[0,+\infty)$. Then for every $x,x'\in X$, $d(\delta_x,\delta_{x'})=\r_{L_d}(\delta_x,\delta_{x'})$.
\end{lemma}
(We remark that Lemma \ref{l5} is different from part i) of Proposition \ref{p1}.)
\begin{proof}
Let $x,x'$ be in $X$. By the definition of $\r_{L_d}$, we have
\begin{equation}\label{e3}
\r_{L_d}(\delta_x,\delta_{x'})=\sup\{|a(x)-a(x')|:\hspace{2mm}a\in\mC,L_d(a)\leq1\}.
\end{equation}
Let $a\in\mC$ and $L_d(a)\leq1$. If $d(\delta_x,\delta_{x'})=0$, then $a(x)=a(x')$ since the map
$\delta_x\longmapsto\delta_x(a)=a(x)$ is continuous with $d$, thus (\ref{e3}) implies that
$$\r_{L_d}(\delta_x,\delta_{x'})=d(\delta_x,\delta_{x'})=0.$$
Now, suppose that $d(\delta_x,\delta_{x'})\neq0$. Since $1=L_d(a)\geq\frac{|a(x)-a(x')|}{d(\delta_x,\delta_{x'})}$,
we have $d(\delta_x,\delta_{x'})\geq|a(x)-a(x')|$, thus (\ref{e3}) implies that
$\r_{L_d}(\delta_x,\delta_{x'})\leq d(\delta_x,\delta_{x'})$. Now, define a map $b_x$ on $X$ by
$b_x(y)=d(\delta_x,\delta_{y})$. Then $b_x\in\mC$ and $L_d(b_x)\leq1$. Thus
$$\r_{L_d}(\delta_x,\delta_{x'})\geq|b_x(x)-b_x(x')|=d(\delta_x,\delta_{x'}).$$
This completes the proof.
\end{proof}
\begin{theorem}
Let $(X,\tau)$, $(Y,\tau')$, $(Z,\tau'')$ be compact Hausdorff spaces and let $d_0$ be a pseudo-metric on $X$ such that
$\tau_{d_0}\subset\tau$. Let
$$F:Y\times Z\to X$$
be a continuous map with $\tau,\tau',\tau''$, and define a pseudo-metric $d_1$ on $Z$ by
$$d_1(z,z')=\sup_{y\in Y}d_0(F(y,z),F(y,z')).$$
With the canonical identification $\bC(Y\times Z,\tau'\times\tau'')\cong\bC(Y,\tau')\ot\bC(Z,\tau'')$ let
$$\hat F:\bC(X,\tau)\to\bC(Y,\tau')\ot\bC(Z,\tau'')$$
be defined by $\hat F(a)=a F$, for $a\in\bC(X,\tau)$. Let
$$(\bC(Z,\tau''),\mC,N)$$
be the QSM space induced by QSM space $(\bC(X,\tau),{\bf Lip}(X,d_0),\|\cdot\|_{d_0})$ and quantum family of
morphisms $(\bC(Z,\tau''),\hat F)$.
Then the following are satisfied.
\begin{enumerate}
\item [i)] $d_1(z,z')=\r_N(\delta_z,\delta_{z'})$ for every $z,z'\in Z$.
\item [ii)]  $\mC\subset{\bf Lip}(Z,d_1)$.
\item [iii)] $\|\cdot\|_{d_1}\leq N$.
\end{enumerate}
\end{theorem}
\begin{proof}
i) Let $L=\|\cdot\|_{d_0}$. Let us recall the definition of $(\bC(Z,\tau''),\mC,N)$. Let $d$ be the pseudo-metric on
$S(\bC(Z,\tau''))$ defined by
$$d(\n,\n')=\sup\{\r_L((\m\ot\n)\hat F,(\m\ot\n')\hat F):\hspace{3mm}\m\in S(\bC(Y,\tau'))\}.$$
Then $N=L_d$ and $\mC$ is the space of all $c\in\bC(Z,\tau'')$ such that the map $\n\longmapsto\n(c)$ on
$S(\bC(Z,\tau''))$ is continuous with $d$ and $N(c)<\infty$. By Lemma \ref{l5}, we have,
\begin{equation}\label{e4}
d(\delta_z,\delta_{z'})=\r_N(\delta_z,\delta_{z'}),
\end{equation}
for every $z,z'\in Z$. Now, we explain the relation between $d_1$ and $d$.\\
Let $z,z'\in Z$ and $y\in Y$. Then
$$(\delta_y\ot\delta_{z})\hat F=\delta_{F(y,z)}\hspace{5mm}\text{and}\hspace{5mm}(\delta_y\ot\delta_{z'})\hat
F=\delta_{F(y,z')}.$$
On the other hand, by Proposition \ref{p1}, for every $x,x'\in X$, $d_0(x,x')=\r_L(\delta_x,\delta_{x'})$.
Thus
$$\r_L((\delta_y\ot\delta_{z})\hat F,(\delta_y\ot\delta_{z'})\hat F)=d_0(F(y,z),F(y,z')).$$
This formula together with the definitions of $d$ and $d_1$, show that
\begin{equation}\label{e5}
d_1(z,z')\leq d(\delta_z,\delta_{z'}).
\end{equation}
Let $\m\in S(\bC(Y,\tau'))$ be arbitrary. We consider $\m$ as a probability Borel regular measure on $(Y,\tau')$.
Then for every $a\in{\bf Lip}(X,d_0)$ with $\|a\|_{d_0}\leq1$, we have,
\begin{equation}\label{e6}
\begin{split}
|(\m\ot\delta_{z})\hat F(a)-(\m\ot\delta_{z'})\hat F(a)|&=|\int_Y(a F(y,z)-a F(y,z'))d_\m(y)|\\
&\leq\int_Y|a(F(y,z))-a(F(y,z'))|d_\m(y).
\end{split}
\end{equation}
For every $y\in Y$, by Lemma \ref{l1},
$$|a(F(y,z))-a(F(y,z'))|\leq d_0(F(y,z),F(y,z')).$$
Therefore, we have
\begin{equation}\label{e7}
|a(F(y,z))-a(F(y,z'))|\leq d_1(z,z').
\end{equation}
(\ref{e7}) and (\ref{e6}) implies that
$$|(\m\ot\delta_{z})\hat F(a)-(\m\ot\delta_{z'})\hat F(a)|\leq d_1(z,z').$$
Therefore, by the definition of $d$,
\begin{equation}\label{e8}
d(\delta_z,\delta_{z'})\leq d_1(z,z').
\end{equation}
Now, by (\ref{e8}) and (\ref{e5}), $d(\delta_z,\delta_{z'})=d_1(z,z')$, and thus by (\ref{e4}),
$$d_1(\delta_z,\delta_{z'})=\r_N(\delta_z,\delta_{z'})$$
for every $z,z'\in Z$, and i) is satisfied.
ii) and iii) are immediate consequence of i) and definitions of $\mC$, $\|\cdot\|_{d_1}$ and $N$.
\end{proof}
\bibliographystyle{amsplain}

\end{document}